\documentclass[reqno,12pt]{amsart}
\usepackage{amsmath}

\newcommand{\CC}{\ensuremath{\mathbb{C}}}

\newcommand{\hgot}{\ensuremath{\mathfrak{h}}}
\newcommand{\lgot}{\ensuremath{\mathfrak{l}}}
\newcommand{\ggot}{\ensuremath{\mathfrak{g}}}
\newcommand{\mgot}{\ensuremath{\mathfrak{m}}}

\newcommand{\sgot}{\ensuremath{\mathfrak{s}}}
\newcommand{\tgot}{\ensuremath{\mathfrak{t}}}

\newtheorem{thm}{Theorem}

\newtheorem{lem}{Lemma}

\theoremstyle{definition}

\theoremstyle{remark}

\begin{document}

\title[Levi-Civita connections of flag manifolds]{Levi-Civita connections of flag manifolds}

\author[A. Sakovich]{Anna Sakovich}

\address{Anna Sakovich \newline Faculty of Pre-University Preparation \newline
Belarusian State University \newline Oktyabrskaya Str.~4 \newline
Minsk 220030, BELARUS } \email{anya\_sakovich@tut.by}

\subjclass[2000]{53C30, 53B05, 17B20}

\keywords{Flag manifold, Levi-Civita connection, invariant
Riemannian metric, root space.}

\thanks{The author will be very grateful for any comments on this paper,
 especially for those concerning the originality of the result obtained.}

\begin{abstract}
For any flag manifold $G/T$ we obtain an explicit expression of
its Levi-Civita connection with respect to any invariant
Riemannian metric.
\end{abstract}
\maketitle

\section{Introduction}

Let $G/T$ be a flag manifold, where $T$ is a maximal torus of a
compact semi-simple Lie group $G$. In this case we obtain an
explicit formula of its Levi-Civita connection (in terms of the
root decomposition for the Lie algebra $\ggot$ of $G$) with
respect to any invariant Riemannian metric. It is possible to
realize this formula, for example, in the case of any classical
simple Lie group $G$. In this paper it is done for $SU(n)$.

This result may prove useful in solving different problems. For
instance, it enables us to determine whether a given metric
$f$-structure $(f,g)$ on $G/T$ belongs to the main classes of
generalized Hermitian geometry (see, for example, \cite{Ki} and
\cite{B}).

\section{Levi-Civita connections of flag manifolds}

In this paper we consider a flag manifold $G/T$, where $T$ is a
maximal torus of a compact semi-simple Lie group $G$. Let $\ggot$
and $\tgot$ be the corresponding Lie algebras of $G$ and $T$.
$G/T$ is a reductive homogeneous space, its reductive
decomposition being $\ggot=\tgot \oplus \mgot$, where $\mgot$ is
an orthogonal complement of $\tgot$ in $\ggot$ with respect to the
Killing form $B$ of $\ggot$. Denote by $\ggot ^\CC$ and $\tgot
^\CC$  the complexifications of $\ggot$ and $\tgot$. Then $\tgot
^\CC$ is a Cartan subalgebra of $\ggot ^\CC$ and we denote by $R$
the root system of $\ggot ^\CC$ with respect to $\tgot ^\CC$. In
this way we have the root decomposition

\begin{equation}\label{rootdecomp} \ggot
^{\CC}=\tgot^{\CC} \oplus \sum_{\alpha \in R} \ggot ^{\alpha}.
\end{equation}

Let $\Pi=\{\alpha_1,\alpha_2,\dots,\alpha_l\}$ be a basis of $R$.
Denote by $R^+$ the set of all positive roots and by $R^-$ the set
of all negative roots. In this paper the following notation will
be used:

\begin{equation*}
|\alpha|=\left\{\begin{array}{ll} \alpha, \; &\text{if} \; \alpha
\in R ^+,
\\
-\alpha,\; & \text{if} \; \alpha \in R ^-.
\end{array}\right.
\end{equation*}

Recall that we can consider the lexicographic order on $R$:
$\gamma=\sum _{i=1} ^n \gamma _i \alpha _i$ is said to be greater
than $\delta=\sum _{i=1} ^n \delta _i \alpha _i$ ($\gamma>\delta$)
if the first nonzero coefficient  $\gamma _k-\delta _k$ in the
decomposition $\gamma-\delta=\sum _{i=1} ^n ( \gamma _i-\delta _i)
\alpha _i$ is positive. If $\gamma-\delta \in R$ then
$\gamma>\delta$ if and only if $\gamma-\delta\in R^+$.

It is well-known that in the case under consideration the
reductive complement $\mgot$ can be decomposed into the direct sum
of 2-dimensional $Ad(T)$-modules $\mgot ^ \alpha$ which are
mutually non-equivalent:

\begin{equation*} \mgot=\sum _{\alpha \in R^+} \mgot^{\alpha},
\text{ where } \mgot ^{\alpha}=\ggot ^{\alpha} \oplus \ggot
^{-\alpha}.
\end{equation*}
Therefore, any invariant Riemannian metric $g=\langle \cdot ,
\cdot \rangle$ on $G/T$ is given by

\begin{equation}\label{metric}
g=\langle \cdot , \cdot \rangle=\sum _{\alpha \in R ^{+}}c
_{\alpha} (\cdot , \cdot)\mid _{\ggot ^{\alpha}\oplus \ggot
^{-\alpha}},
\end{equation}
where $c_{\alpha}>0$, $(\cdot , \cdot)$ is the negative of the
Killing form $B$ of the Lie algebra $\ggot$.

In this paper we will need the following result.

\begin{thm}\cite{KN} Let $(M,g)$ be a Riemannian manifold, $M=G/H$
a reductive homogeneous space with the reductive decomposition
$\ggot=\hgot \oplus \mgot$. Then the Levi-Civita connection with
respect to $g$ can be expressed in the form
\begin{equation}\label{LC}
\nabla_X Y=\frac{1}{2}[X,Y]_{\mgot}+U(X,Y),
\end{equation}
where $U$ is a symmetric bilinear mapping $\mgot \times \mgot
\rightarrow \mgot$ defined by the formula
\begin{equation}\label{U}
2g(U(X,Y),Z)=g(X,[Z,Y]_{\mgot})+g([Z,X]_{\mgot},Y), \; X,Y,Z \in
\mgot.
\end{equation}
\end{thm}

We can consider (\ref{U}) as an equation of variable $U$. Let us
try to solve this equation in the case of an arbtrary flag
manifold $G/T$.

We begin with obtaining an important preliminary result. Consider
$X_{\gamma} \in \ggot ^{\gamma}$, $Y_{\delta} \in \ggot
^{\delta}$, $\gamma$, $\delta \in R$.  In the view of
(\ref{metric}), (\ref{U}) takes the following form:

\begin{multline}\label{!}
2 \sum _{\alpha \in R ^+} c_{\alpha}
(U(X_{\gamma},Y_{\delta})_{\ggot^{\alpha}\oplus
\ggot^{-\alpha}},Z_{\ggot^{\alpha} \oplus \ggot^{-\alpha}})\\=
\sum_{\alpha \in R ^+} c_{\alpha}
((X_{\gamma})_{\ggot^{\alpha}\oplus
\ggot^{-\alpha}},[Z,Y_{\delta}]_{\ggot^{\alpha}\oplus
\ggot^{-\alpha}})\\+\sum_{\alpha \in R ^+} c_{\alpha}
([Z,X_{\gamma}]_{\ggot^{\alpha}\oplus
\ggot^{-\alpha}},(Y_{\delta})_{\ggot^{\alpha}\oplus
\ggot^{-\alpha}}).
\end{multline}
Obviously, the right-hand side of this equation is equal to
\begin{equation*}
c_{|\gamma|}(X_{\gamma},[Z,Y_{\delta}]_{\ggot^{|\gamma|}\oplus\ggot^{-|\gamma|}})+
c_{|\delta|}([Z,X_{\gamma}]_{\ggot^{|\delta|}\oplus\ggot^{-|\delta|}},Y_{\delta}).
\end{equation*}

Let $Z=\sum _{\alpha \in R}Z_{\alpha}$, where $Z_{\alpha}=Z_{\ggot
^{\alpha}}$. Note that
\begin{equation*}
[Z,Y_{\delta}]_{\mgot}=[\sum _{\alpha \in R}
Z_{\alpha},Y_{\delta}]_{\mgot}=\sum _{\alpha \in R}
[Z_{\alpha},Y_{\delta}]_{\mgot}=\sum _{\alpha, \alpha+\delta \in
R} [Z_{\alpha},Y_{\delta}],
\end{equation*}
and, evidently, $[Z_{\alpha},Y_{\delta}]=[Z,Y_{\delta}]_{\ggot
^{\alpha \oplus \delta}}$. It is easy to see that
\begin{multline*}
(X_{\gamma},[Z,Y_{\delta}]_{\ggot^{|\gamma|}\oplus\ggot^{-|\gamma|}})=\left
(X_{\gamma},\left( \sum_{\alpha,\alpha+\delta\in
R}[Z_{\alpha},Y_{\delta}]\right)
_{\ggot^{|\gamma|}\oplus\ggot^{-|\gamma|}}\right)\\=
\left(X_{\gamma},\left( \sum_{\alpha,\alpha+\delta\in
R}[Z_{\alpha},Y_{\delta}]\right) _{\ggot^{-\gamma}}\right).
\end{multline*}

If $ \left( \sum_{\alpha,\alpha+\delta\in
R}[Z_{\alpha},Y_{\delta}]\right) _{\ggot^{-\gamma}}\neq 0$, then
there exists such $\alpha \in R$ that $\alpha+\delta=-\gamma$. In
other words, $\alpha=-\gamma-\delta \in R$. Therefore,
\begin{equation*}
c_{|\gamma|}(X_{\gamma},[Z,Y_{\delta}]_{\ggot^{|\gamma|}\oplus\ggot^{-|\gamma|}})=
\left\{\begin{array}{ll} 0, \; &\text{if } \; \gamma+\delta \notin
R,\\ c_{|\gamma|}(X_{\gamma},[Z_{-\gamma-\delta},Y_{\delta}]),\; &
\text{if } \;  \gamma+\delta \in R.
 \end{array}\right.
\end{equation*}

Arguing as above, one can prove that
\begin{equation*}
c_{|\delta|}([Z,X_{\gamma}]_{\ggot^{|\delta|}\oplus\ggot^{-|\delta|}},Y_{\delta})=
\left\{\begin{array}{ll} 0, \; &\text{if } \; \gamma+\delta \notin
R,\\ c_{|\delta|}([Z_{-\delta-\gamma},X_{\gamma}],Y_{\delta}),\; &
\text{if } \; \gamma+\delta \in R.
 \end{array}\right.
\end{equation*}

Hence, if $\gamma+\delta \notin R$, (\ref{U}) is transformed into
\begin{equation*}
2g(U(X_{\gamma},Y_{\delta}),Z)=0
\end{equation*}
for any $Z \in \mgot$. Thus, if $\gamma+\delta \notin R$, then
$U(X_{\gamma},Y_{\delta})=0$.

If $\gamma+\delta \in R$, then (\ref{!}) is equivalent to
\begin{multline*}
 2\sum _{\alpha \in R ^+} c_{\alpha}
(U(X_{\gamma},Y_{\delta})_{\ggot^{\alpha}\oplus
\ggot^{-\alpha}},Z_{\ggot^{\alpha} \oplus
\ggot^{-\alpha}})\\=c_{|\gamma|}(X_{\gamma},[Z_{-\gamma-\delta},Y_{\delta}])
+c_{|\delta|}([Z_{-\delta-\gamma},X_{\gamma}],Y_{\delta}).
\end{multline*}

By the properties of the Killing form we obtain
\begin{multline*}
 \sum _{\alpha \in R ^+ \atop \alpha \neq|\gamma+\delta|} (2c_{\alpha}
U(X_{\gamma},Y_{\delta})_{\ggot^{\alpha}\oplus
\ggot^{-\alpha}},Z_{\ggot^{\alpha} \oplus \ggot^{-\alpha}})+
(2c_{|\gamma+\delta|}U(X_{\gamma},Y_{\delta})_{\ggot^{|\gamma+\delta|}
\oplus
\ggot^{-|\gamma+\delta|}}\\-c_{|\gamma|}[Y_{\delta},X_{\gamma}]-c_{|\delta|}
[X_{\gamma},Y_{\delta}],Z_{\ggot^{|\gamma+\delta|} \oplus
\ggot^{-|\gamma+\delta|}})=0.
\end{multline*}

Since $\mgot^{\alpha}$ is orthogonal to $\mgot^{\beta}$  with
respect to the Killing form of $\ggot$ ($\alpha$, $\beta \in R^+$,
$\alpha \neq \beta$),  we have
\begin{multline*}
(2\sum _{\alpha \in R ^+ \atop \alpha
\neq|\gamma+\delta|}c_{\alpha}
U(X_{\gamma},Y_{\delta})_{\mgot^{\alpha}}+2c_{|\gamma+\delta|}U(X_{\gamma},Y_{\delta})_{\mgot^{|\gamma+\delta|}}
\\-(c_{|\gamma|}-c_{|\delta|})[Y_{\delta},X_{\gamma}],Z)
 =0
\end{multline*}
for any $Z \in \mgot$. This yields that $$2\sum _{\alpha \in R ^+
\atop \alpha \neq|\gamma+\delta|}c_{\alpha}
U(X_{\gamma},Y_{\delta})_{\mgot^{\alpha}}+2c_{|\gamma+\delta|}U(X_{\gamma},Y_{\delta})_{\mgot^{|\gamma+\delta|}}
\\-(c_{|\gamma|}-c_{|\delta|})[Y_{\delta},X_{\gamma}]$$
(and, consequently, any of its projections onto $\mgot^{\alpha}$,
$\alpha \in R^+$) is equal to 0. We have proved the following
result.

\begin{lem}
Let $G/T$ be a flag manifold with the root decomposition
(\ref{rootdecomp}). Then for any $X_{\gamma} \in \ggot ^{\gamma}$,
$Y _{\delta} \in \ggot ^{\delta}$, where $\gamma$, $\delta \in R$,
we have
\begin{equation}\label{U1}
U(X_{\gamma},Y_{\delta})=\left\{\begin{array}{ll}
\frac{c_{|\gamma|}-c_{|\delta|}}{2c_{|\gamma+\delta|}}[Y_{\delta},X_{\gamma}],
\; &\text{ if } \; \gamma + \delta \in R,
\\
0,\; & \text{ if } \; \gamma + \delta \notin R.
\end{array}\right.
\end{equation}
\end{lem}

This lemma enables us to obtain the similar expression for
$U(X,Y)$ in the case of any $X=\sum _{\alpha \in R} X_{\alpha}$
and $Y=\sum _{\beta \in R} Y_{\beta}$ in $\mgot$. As $U$ is
bilinear, application of (\ref{U1}) gives us
\begin{equation}\label{a}
U(X,Y)=\sum _{\alpha, \beta \in R} U(X_{\alpha},Y_{\beta})=\sum
_{\alpha, \beta , \alpha+\beta \in R}
\frac{c_{|\alpha|}-c_{|\beta|}}{2c_{|\alpha+\beta|}}[Y_{\beta},X_{\alpha}].
\end{equation}

For any $\alpha, \beta \in R$ such that $\alpha+\beta \in R$ we
group together terms with the coefficient
$\frac{c_{|\alpha|}-c_{|\beta|}}{2c_{|\alpha+\beta|}}$. In this
way we obtain the sum of the following summands
\begin{equation*}
\frac{c_{|\alpha|}-c_{|\beta|}}{2c_{|\alpha+\beta|}}Z_{\alpha}
^{\beta},
\end{equation*}
where
\begin{equation}\label{Z}
Z_{\alpha} ^{\beta}=[Y_{\beta},
X_{\alpha}]+[X_{\beta},Y_{\alpha}]+[Y_{-\beta},
X_{-\alpha}]+[X_{-\beta},Y_{-\alpha}],\; \alpha, \beta \in R.
\end{equation}

However, $Z_{\alpha} ^{\beta}=Z_{-\alpha} ^{-\beta}=Z_{\beta}
^{\alpha}=Z_{-\beta} ^{-\alpha}$, which implies that there is a
need to restrict the range of $\alpha$ and $\beta$. Certainly,
(\ref{a}) is equivalent to
\begin{equation*}
U(X,Y)=\frac{1}{4} \sum _{\alpha, \beta \in R}
\frac{c_{|\alpha|}-c_{|\beta|}}{2c_{|\alpha+\beta|}} Z_{\alpha}
^{\beta},
\end{equation*}
but this formula is definitely not the most convenient since there
are repetitions of summands. Luckily, it is easy to establish a
condition which makes it possible to select one pair of roots out
of four pairs $(\alpha,\beta)$, $(\beta,\alpha)$,
$(-\alpha,-\beta)$, $(-\beta,-\alpha)$.

\begin{lem}
For any $\alpha, \beta \in R$ there exists only one pair
$(a_1,a_2)
\in\{(\alpha,\beta),(\beta,\alpha),(-\alpha,-\beta),(-\beta,-\alpha)\}$
such that $|a_1|<a_2$.
\end{lem}

\begin{proof}
The condition $|a_1|<a_2$ presupposes that $a_2 \in R ^+$.
Obviously, $|a_1|<a_2$ if and only if $-a_2<a_1<a_2$.

Such a pair can be chosen as follows.

Set $a_1=\alpha$, $a_2=\beta$. If $a_2\in R^-$, set $a_1$ equal to
$-a_1$ and $a_2$ equal to $-a_2$. Thus we have $a_2\in R^+$. Now
let us check if $a_1<a_2$. If this condition is not satisfied, set
$a_2$ equal to $a_1$ and $a_1$ equal to $a_2$. It remains to
verify if $a_1>-a_2$. If this is true, the desired pair
$(a_1,a_2)$ is obtained , otherwise we choose $(-a_2,-a_1)$.

The uniqueness of this pair can be proved as follows. Without loss
of generality, suppose that $|\alpha|<\beta$, that is,
$-\beta<\alpha<\beta$. Then $(\beta,\alpha)$ satisfies
$\beta>\alpha$ and for $(-\alpha,-\beta)$ we have $-\alpha>-\beta$
which means that these two pairs do not satisfy the stipulated
condition. The pair $(-\beta,-\alpha)$ should satisfy
$\alpha<-\beta<-\alpha$ and this contradicts the assumption made
above.
\end{proof}

\medskip
In the view of this lemma we have
\begin{equation}\label{2}
U(X,Y)=\sum _{\alpha,\beta,\alpha+\beta \in R,
 \atop |\alpha|<\beta \in R
 ^+}\frac{c_{|\alpha|}-c_{|\beta|}}{2c_{|\alpha+\beta|}}Z_{\alpha}^{\beta}
\end{equation}
($Z_{\alpha}^{\beta}$ is determined by means of (\ref{Z})).

Consider different cases for $\alpha, \beta \in R$. $\beta$ always
belongs to $R ^+$ and $\alpha$  can be selected from both $R^+$
and $R^-$.

If $\alpha \in R ^+$, $\beta \in R^+$ then the conditions
$\alpha+\beta \in R$ and $|\alpha|<\beta$ can be replaced by the
conditions $\alpha+\beta \in R ^+$ and $\alpha<\beta$
respectively.

If $\alpha \in R ^-$, $\beta \in R^+$ then $|\alpha|<\beta$ is
equivalent to $-\alpha<\beta$. If $\alpha+\beta \in R$ then
$-\alpha<\beta$ can be substituted for the condition $\alpha+\beta
\in R ^+$.

Therefore, the right-hand side of (\ref{2}) is transformed into
\begin{multline*}
\sum _{\alpha,\beta,\alpha+\beta \in R ^+,
 \atop
 \alpha<\beta}\frac{c_{\alpha}-c_{\beta}}{2c_{\alpha+\beta}}Z_{\alpha}^{\beta}+
 \sum _{-\alpha,\beta,\alpha+\beta \in R
 ^+}\frac{c_{-\alpha}-c_{\beta}}{2c_{\alpha+\beta}}Z_{-\alpha}^{\beta}\\=\sum _{\alpha,\beta,\alpha+\beta \in R ^+,
 \atop
 \alpha<\beta}\frac{c_{\alpha}-c_{\beta}}{2c_{\alpha+\beta}}Z_{\alpha}^{\beta}+
 \sum _{\alpha,\beta,\beta-\alpha \in R
 ^+}\frac{c_{\alpha}-c_{\beta}}{2c_{\beta-\alpha}}Z_{\alpha}^{\beta}.
\end{multline*}

Thus, the following theorem is proved.

\begin{thm} Let $G/T$ be a flag manifold with the root decomposition $(\ref{rootdecomp})$.
Then for any $X$, $Y\in \mgot$ we have
\begin{equation}\label{3}
U(X,Y)=\sum _{\alpha,\beta,\alpha+\beta \in R ^+,
 \atop
 \alpha<\beta}\frac{c_{\alpha}-c_{\beta}}{2c_{\alpha+\beta}}Z_{\alpha}^{\beta}+
 \sum _{\alpha,\beta,\beta-\alpha \in R
 ^+}\frac{c_{\alpha}-c_{\beta}}{2c_{\beta-\alpha}}Z_{\alpha}^{\beta},
\end{equation}
where $Z_{\alpha} ^{\beta}=[Y_{\beta},
X_{\alpha}]+[X_{\beta},Y_{\alpha}]+[Y_{-\beta},
X_{-\alpha}]+[X_{-\beta},Y_{-\alpha}],\; \alpha, \beta \in R.$
\end{thm}

\section{Examples}

As an example, let us consider the flag manifold $G/T=SU(n+1)/T$
$(n \geq 2 )$, where $T$ is a maximal torus of $SU(n+1)$.

In this case $\ggot=\sgot \lgot (n+1, \CC)$. The root system of
$SU(n+1)$ with respect to $\tgot$ is
\begin{equation*}
R=A_n=\{\varepsilon _i-\varepsilon _j \; | \; i \neq j, \; 1\leq
i,j \leq n+1\},
\end{equation*}
its basis being
\begin{equation*}
\{\alpha _i=\varepsilon _i-\varepsilon _{i+1}\}_{1\leq i \leq n}.
\end{equation*}
The set of all positive roots in this case is
\begin{equation*}
R ^+=\{\varepsilon _i-\varepsilon _j \; | \; 1 \leq i<j \leq n\}.
\end{equation*}

An arbitrary positive root $\alpha=\varepsilon _i-\varepsilon _j$,
where $i<j$, is decomposed into the sum of basis vectors as
follows:
\begin{equation*}
\alpha=\varepsilon _i-\varepsilon _j=\alpha _i+\alpha_{i+1}+ \dots
+\alpha _j.
\end{equation*}
It is easy to see that $\alpha=\varepsilon _i-\varepsilon
_j<\beta=\varepsilon _k-\varepsilon _l$ ($\alpha,\beta \in R ^+$)
if and only if $i>k$.

Take $\alpha=\varepsilon _i-\varepsilon _j$, $\beta=\varepsilon
_k-\varepsilon _l \in R ^+$, where $i<j$, $k<l$.

$\alpha+\beta \in R ^+$ if and only if either $i<j=k<l$ (hence
$\alpha+\beta=\varepsilon _i-\varepsilon _l$), or $k<i=l<j$ (hence
$\alpha+\beta=\varepsilon _k-\varepsilon _j$). Note that in the
first case $\alpha>\beta$ and in the second case $\beta>\alpha$.

$\beta-\alpha \in R ^+$ if and only if either $i=k<j<l$ (hence
$\beta-\alpha=\varepsilon _j-\varepsilon _l$), or $k<i<j=l$ (hence
$\beta-\alpha=\varepsilon _k-\varepsilon _i$).

It is not difficult to show that $Z_{\alpha} ^{\beta}=[X_{\mgot
^{\beta}},Y_{\mgot ^{\alpha}}]+[Y_{\mgot ^{\beta}},X_{\mgot
^{\alpha}}]$ for any $\alpha, \beta \in R^+$.

Therefore, in the case of $SU(n+1)/T_{max}$ $(n \geq 2)$ (\ref{3})
takes form
\begin{multline}\label{SU}
U(X,Y)\\=\sum _{1 \leq i <j <k \leq
n+1}\frac{c_{\varepsilon_j-\varepsilon_k}-c_{\varepsilon_i-\varepsilon_j}}
{2c_{\varepsilon_i-\varepsilon_k}}([X_{\mgot
^{\varepsilon_i-\varepsilon_j}},Y_{\mgot
^{\varepsilon_j-\varepsilon_k}}]+[Y_{\mgot ^{
\varepsilon_i-\varepsilon_j}},X_{\mgot
^{\varepsilon_j-\varepsilon_k}}])
\\
 +\sum _{1 \leq i <j <k \leq
n+1}\frac{c_{\varepsilon_i-\varepsilon_j}-c_{\varepsilon_i-\varepsilon_k}}
{2c_{\varepsilon_j-\varepsilon_k}}([X_{\mgot
^{\varepsilon_i-\varepsilon_k}},Y_{\mgot
^{\varepsilon_i-\varepsilon_j}}]+[Y_{\mgot
^{\varepsilon_i-\varepsilon_k}},X_{\mgot
^{\varepsilon_i-\varepsilon_j}}])\\+ \sum _{1 \leq i <j <k \leq
n+1}\frac{c_{\varepsilon_j-\varepsilon_k}-c_{\varepsilon_i-\varepsilon_k}}
{2c_{\varepsilon_i-\varepsilon_j}} ([X_{\mgot
^{\varepsilon_i-\varepsilon_k}},Y_{\mgot
^{\varepsilon_j-\varepsilon_k}}]+[Y_{\mgot
^{\varepsilon_i-\varepsilon_k}},X_{\mgot
^{\varepsilon_j-\varepsilon_k}}]).
\end{multline}

As a particular case, let us consider the flag manifold
$SU(3)/T_{max}$. The set of all positive roots is
\begin{equation*}
R ^+=\{\alpha_1=\varepsilon _1-\varepsilon _2,\;
\alpha_2=\varepsilon _1-\varepsilon _3,\; \alpha_3=\varepsilon
_2-\varepsilon _3 \}.
\end{equation*}

In order to obtain a more compact formula denote $c_{\alpha _i}$
by $c_i$ and $\mgot ^{\alpha _i}$ by $\mgot _i$. We also agree to
write $X_i$ instead of $X_{\mgot _i}$.

Therefore, in the case  of $SU(3)/T_{max}$, using the notations
introduced above, we can rewrite (\ref{SU}) as follows:
\begin{multline*}
U(X,Y)=\frac{c_3-c_2}{2c_1}([X_2,Y_3]+[Y_2,X_3])\\
+\frac{c_3-c_1}{2c_2}([X_1,Y_3]+[Y_1,X_3])+\frac{c_2-c_1}{2c_3}([X_1,Y_2]+[Y_1,X_2]).
\end{multline*}
Actually, this result is well-known (see, for example, \cite{W}).

\newpage
\renewcommand{\refname}{REFERENCES }
\begin{center}

\end{center}
\end{document}